\documentclass{aims}
\usepackage{amsmath}
  \usepackage{paralist}
  \usepackage{graphics} 
  \usepackage{epsfig} 
 \usepackage[colorlinks=true]{hyperref}
\hypersetup{urlcolor=blue, citecolor=red}

\def\currentvolume{X}
 \def\currentissue{X}
  \def\currentyear{200X}
   \def\currentmonth{XX}
    \def\ppages{X--XX}

\newtheorem{theorem}{Theorem}[section]

\newtheorem{lemma}[theorem]{Lemma}
\newtheorem{proposition}{Proposition}

\theoremstyle{definition}
\newtheorem{definition}[theorem]{Definition}
\newtheorem{remark}{Remark}

\newcommand{\field}[1]{\mathbb{#1}}

\newcommand{\RR}{\field{R}}

\def\id{\mathop{\rm id}\nolimits}

\def\vect#1{\overrightarrow{#1}}
\title[Conformally Hamiltonian vector fields]
      {A property of conformally Hamiltonian vector fields; application to the Kepler problem}

\author[Charles-Michel Marle]{}


 \email{marle@math.jussieu.fr,\quad charles-michel.marle@polytechnique.org}
\thanks{}

\begin{document}
\maketitle

\centerline{\scshape Charles-Michel Marle }
\medskip
{\footnotesize
 \centerline{Universit\'e Pierre et Marie Curie}
 \centerline{Institut de math\'ematiques de Jussieu}
   \centerline{4 place Jussieu, case courrier 247}
   \centerline{75252 Paris cedex 05, France}
} 

\bigskip

 \centerline{(Communicated by the associate editor name)}

\begin{abstract}
Let $X$ be a Hamiltonian vector field defined on a symplectic manifold $(M,\omega)$, $g$ a nowhere vanishing
smooth function defined on an open dense subset $M^0$ of $M$. We will say that the vector field $Y=gX$ is
\emph{conformally Hamiltonian}. We prove that when $X$ is complete, when $Y$ is Hamiltonian with respect to
another symplectic form $\omega_2$ defined on $M^0$, and when another technical condition is satisfied,  
then there is a symplectic diffeomorphism from $(M^0,\omega_2)$ onto an open subset of $(M,\omega)$, 
which maps each orbit to itself and is equivariant with respect to the flows of the vector fields 
$Y$ on $M^0$ and $X$ on $M$. This result explains why the diffeomorphism of the phase space of the Kepler problem restricted 
to the negative (resp. positive) values of the energy function, onto an open subset of the cotangent bundle
to a three-dimensional sphere (resp. two-sheeted hyperboloid), discovered by Gy\"orgyi (1968)~\cite{gyorgyi}, 
re-discovered by  Ligon and Schaaf (1976)~\cite{ligonschaaf}, is a symplectic diffeomorphism.  
Cushman and Duistermaat (1997)~\cite{cushmanduistermaat} have shown that the Gy\"orgyi-Ligon-Schaaf diffeomorphism is characterized
by three very natural properties; here that diffeomorphism is obtained by composition of the diffeomorphism
given by our result about conformally Hamiltonian vector fields with a (non-symplectic) 
diffeomorphism built by a variant of Moser's method~\cite{moser}. Infinitesimal 
symmetries of the Kepler problem are discussed, and it is shown that their space is a Lie algebroid
with zero anchor map rather than a Lie algebra.
\end{abstract}

\section{Introduction}
I am very glad to submit a paper for the special issue of the Journal of Geometric Mechanics in honour of
Tudor Ratiu. I followed his scientific work for several years; I specially praise the superb book he wrote
with Juan-Pablo Ortega~\cite{ortegaratiu}. I hope that he will find some interest in the present work.

The Kepler problem is a completely integrable Hamiltonian system with important applications 
in the physical world: it provides a very accurate model of the motion of planets in the solar system, 
and its quantized equivalent provides a good model of the hydrogen atom.
Moreover, several features of the Kepler problem make it very interesting for the mathematician:
some of its motions (those in which there is a collision of the moving point with the attractive centre) 
are not defined for all times, but the system can be regularized, \emph{i.e.} mapped into a new 
Hamitonian system whose motions are defined for all times;
the Lie algebra of infinitesimal symmetries of a given energy level of the phase space depends on 
that energy level.
\par\smallskip

In Section~\ref{conformallyfields} some results in symplectic geometry are presented. 
These results will be used in Section~\ref{kepler} to explain why a diffeomorphism 
between the phase space of the Kepler problem, restricted to negative values 
(resp. positive values) of the energy, and an open subset of the cotangent space 
to a 3-dimensional sphere (resp., a 3-dimensional two-sheeted hyperboloid)
is symplectic. That diffeomorphism was discovered by Gy\"orgyi~\cite{gyorgyi} (1968), 
re-discovered by Ligon and Schaaf~\cite{ligonschaaf} (1976) and discussed by Cushman and 
Duistermaat~\cite{cushmanduistermaat} (1997), who have shown that it is characterized by the three very natural properties:

\begin{enumerate}

\item It maps the set of points in the phase space of the Kepler problem where the energy is negative (resp., positive)
onto the tangent bundle of the $3$-sphere (resp., the two-sheeted $3$-dimensional revolution hyperboloid) with its zero section removed.

\item It intertwines the Kepler and Delaunay vector fields (a rescaling of the geodesic vector field on the $3$-sphere, 
or on the $3$-dimensional hyperboloid).

\item It intertwines the $so(4)$-momentum mappings of the Kepler and Delauney vector fields.

\end{enumerate}

We will see that the remarkable properties of that
diffeomorphism appear as very natural consequences of the results presented in Section~\ref{conformallyfields}.
\par\smallskip

We will also discuss the weak infinitesimal 
symmetries of the Kepler problem, and we will show that their set is a Lie algebroid 
with a zero anchor map, rather than a Lie algebra.
\par\smallskip

After completion and submission of the present paper, 
the work of Gert~Heckman and Tim~de~Laat~\cite{heckmandelaat}, recently posted on arXiv,
was indicated to us. The method used by these authors to explain the properties of the
Gy\"orgyi-Ligon-Schaaf diffeomorphism rests on the same ideas as ours. We also learnt that conformally Hamiltonian vector fields were used,
in the theory of bi-Hamiltonian vector fields, by A.J.~Maciejewski, M. Prybylska and A.V.~Tsiganov~\cite{MaciePryTsi}.

\section{Conformally Hamiltonian vector fields}\label{conformallyfields}

\subsection{Notations and conventions}\label{notations}
Unless another assumption is explicitly stated, all manifolds, functions, applications, vector fields and differential forms considered in this work will be assumed to be smooth, \emph{i.e} differentiable of class $C^\infty$.

\begin{definition}\label{flow}
Let $X$ be a vector field on a manifold $M$. The \emph{differential equation determined by} $X$ is the ordinary differential equation
\begin{equation*}
 \frac{d\varphi(t)}{dt}=X\bigl(\varphi(t)\bigr)\,.
 \end{equation*}
The \emph{flow} of that equation is     
the map $\Phi_X$, defined on an open subset
$D_X$ of $\RR\times M$, with values in $M$, such that, for each $x\in M$, the parametrized curve 
$t\mapsto \Phi_X(t,x)$ is the maximal solution $\varphi$ of the above differential equation
which takes the value $x$ for $t=0$. It means that
 $$\frac{\partial\Phi(t,x)}{\partial t}=X\bigl(\Phi(t,x)\bigr)
 \quad\text{for each\ }(t,x)\in D_X\,,\quad \Phi_X(0,x)=x\quad\text{for each\ }x\in M\,,$$
and that the open subset $D_X$ of $\RR\times M$ is such that, for each $x\in M$
 $$I_x=\bigl\{t\in\RR\bigm|(t,x)\in D_x\bigr\}$$
is the largest open interval of $\RR$ on which a solution $\varphi$ of the differential equation 
determined by $X$ satisfying $\varphi(0)=x$ can be defined.
\end{definition}

\begin{remark}\label{partial flow} 
The map $\Phi_X$ is sometimes called \emph{partial flow} to distinguish it from the \emph{full flow} $\Psi_X$, which is the map 
defined on an open subset of $\RR\times\RR\times M$, with values in $M$, such that, for each $t_0\in \RR$ and $x \in M$, 
the parametrized curve $t\mapsto \Psi_X(t,t_0,x)$ is the maximal solution of the above differential equation
which takes the value $x_0$ for $t=t_0$. The vector field $X$ being time-independent, we have $\Psi_X(t,t_0 ,x_0)=\Phi_X(t-t_0,x_0)$, 
so the full flow is determined by the partial flow. It is no longer true for ordinary differential equations determined by a 
time-dependent vector field; the full flow should then be used instead of the partial flow.
\end{remark}
 
\subsection{Change of independent variable in ordinary differential equations}\label{changevar}

\begin{lemma}
On a manifold $M$, let $X$ be a vector field and $g$ a nowhere vanishing function. Consider the two 
ordinary differential equations
 \begin{align}
 \frac{d\varphi(t)}{dt}&=g\bigl(\varphi(t)\bigr)X\bigl(\varphi(t)\bigr)\,,\tag{$*$}\\
 \frac{d\psi(s)}{ds}&=X\bigl(\psi(s)\bigr)\,.\tag{$**$}
 \end{align}
Let us assume that there exists a smooth function $\sigma:\RR\times M\to\RR$ 
such that for each solution $\varphi$ of $(*)$
 \begin{equation}
 \frac{d}{dt}\sigma\bigl(t,\varphi(t)\bigr)=g\bigl(\varphi(t)\bigr)\,.\tag{${*}{*}{*}$}
 \end{equation}    
For each solution $\varphi:I_\varphi\to M$ of $(*)$ defined on the open interval $I_\varphi$ of $\RR$, let
 $$\sigma_\varphi(t)=\sigma\bigl(t,\varphi(t))\,.$$
Then the function $\sigma_\phi$ is a diffeomorphism from $I_\varphi$ onto another open interval $\sigma_\varphi(I_\varphi)$
of $\RR$, and the map  
  $$\psi:\sigma_\varphi(I_\varphi)\to M\,,\quad s\mapsto \psi(s)=\varphi\circ\sigma_\varphi^{-1}(s)$$
is a solution of the ordinary differential equation $(**)$. 
\end{lemma}

\begin{proof}
The derivative of the function $\sigma_\varphi$, at each poit $t\in I_\varphi$, 
is $\displaystyle g\bigl(\varphi(t))$, which never vanishes. Therefore $\sigma_\varphi$ is a diffeomorphism.
The chain rule shows that $\psi$ is a solution of $(**)$.
\end{proof}

\begin{remark}
For a general smooth vector field $X$ and a general smooth nowhere vanishing function $g$ given on $M$,
there may be no globally defined function $\sigma:\RR\times M\to\RR$ verifying $({*}{*}{*})$. In the Kepler
problem, that function exists, as we will see in subsection \ref{levi-civitaparameter}, and is affine in 
the variable $t$. The following proposition gives some information about the existence of the map $\sigma$.
\end{remark}

\begin{proposition}
On a manifold $M$, let $X$ be a vector field and $g:M\to\RR\backslash\{0\}$ 
a nowhere vanishing function. We denote by $\Phi_X:D_X\to M$ the flow of $X$~(\ref{flow}).
For each $t_0\in\RR$, there exists a smooth function $\sigma$, defined on an open neighbourhood $W_{t_0}$
of $\{t_0\}\times M$ in $\RR\times M$, such that for each solution $\varphi:I_\varphi\to M$ of the  
differential equation determined by $X$ defined on an open interval $I_\varphi$ containing $t_0$, and each
$t\in I_\varphi$,
\begin{equation}
 \frac{d}{dt}\sigma\bigl(t,\varphi(t)\bigr)=g\bigl(\varphi(t)\bigr)\,.\tag{${*}$}
 \end{equation}    
The function $\sigma$ is not unique: any smooth function defined on $M$ can be chosen for its restriction to
$\{t_0\}\times M$. The function $\sigma$ is affine with respect to the variable $t$ if and only if its  restriction 
$\sigma_{t_0}:M\to \RR$, $\sigma_{t_0}(x)=\sigma(t_0,x)$ is such that the Lie derivative with respect to
the vector field $X$ of the function $g-\sigma_{t_0}$
is constant along each integral curve of $X$. When the preceding condition is satisfied,
the affine extension of $\sigma$ to the whole $\RR\times M$, still denoted by $\sigma$, satisfies $(*)$ for 
all solutions $\varphi:I\to M$ of the differential equation determined by $X$
and all $t\in I$.
\end{proposition}

\begin{proof} 
Let us choose any $t_0\in \RR$.
The map $(\theta,x)\mapsto(t_0-\theta,x)$ is a diffeomorphism of $\RR\times M$ onto itself, which sends
the open subset $D_X$ on which the flow $\Phi_X$ is defined onto the subset 
$W_{t_0}=\bigl\{\,(t,x)\in\RR\times M\bigm|(t_0-t,x)\in D_X\}$. Therefore $W_{t_0}$ 
is an open subset of $\RR\times M$ 
which contains $\{t_0\}\times M$. Let $\sigma_{t_O}:M\to\RR$ be any smooth function. The formula
 $$\sigma(t,x)=\sigma_{t_0}\circ\Phi_X(t_0-t,x)+\int_{t_0}^tg\circ\Phi_X(\tau-t,x)\,d\tau$$
defines a smooth function $\sigma:W_{t_0}\to\RR$ whose restriction to $\{t_0\}\times M$ is $(t_0,x)\mapsto
\bigl(t_0,\sigma_{t_0}(x)\bigr)$. For each $(t,x_0)\in\RR\times M$ such that $(t-t_0,x_0)\in D_X$, we have
 \begin{align*}
 \sigma\bigl(t,\Phi_X(t-t_0,x_0)\bigr)&=\sigma_{t_0}\circ\Phi_X\bigl(t_0-t,\Phi_X(t-t_0,x_0)\bigr)\\
                                      &\quad\quad +\int_{t_0}^tg\circ\Phi_X\bigl(\tau-t,\Phi_X(t-t_0,x_0)\bigr)\,d\tau\\
                                      &=\sigma_{t_0}(x_0)+\int_{t_0}^tg\circ\Phi_X(\tau-t_0,x_0)\,d\tau\,.     
 \end{align*}
Therefore
 $$\frac{d}{dt}\sigma\bigl(t,\Phi(t-t_0,x_0)\bigr)=g\bigl(\Phi_X(t-t_0,x_0)\bigr)\,.$$
Since $t\mapsto\Phi_X(t-t_0,x_0)$ is the maximal solution of the differential equation determined by $X$ which takes the value $x_0$ for $t=t_0$, we see that the map $\sigma$ satisfies condition $(*)$ of the statement above.

The map $\sigma$ is affine with respect to $t$ if and only if its partial derivative with respect to $t$ does not
depend on $t$. We have
 \begin{align*}
 \frac{\partial\sigma(t,x)}{\partial t}
 &=-\bigl\langle d\sigma_{t_0}\circ\Phi_X(t_0-t,x),X\circ\Phi_X(t_0-t,x)\bigr\rangle\\
 &\quad\quad +g(x)+\int_{t_0}^t\frac{\partial}{\partial t}\bigl(g\circ\Phi_X(\tau-t,x)\bigr)\,d\tau\,.
 \end{align*}
Taking into account
 $$\frac{\partial}{\partial t}\bigl(g\circ\Phi_X(\tau-t,x)\bigr)
 =-\frac{\partial}{\partial \tau}\bigl(g\circ\Phi_X(\tau-t,x)\bigr)
 $$
we obtain
 \begin{align*}
 \frac{\partial\sigma(t,x)}{\partial t}
 &=-\bigl\langle d\sigma_{t_0}\circ\Phi_X(t_0-t,x),X\circ\Phi_X(t_0-t,x)\bigr\rangle\\
 &\quad\quad +g(x)-\Bigl(\bigl\langle dg\circ\Phi_X(\tau-t,x),X\circ\Phi_X(\tau-t,x)\bigr\rangle\Bigr)
 \Bigm|_{\tau=t_0}^{\tau=t}\,.
 \end{align*}
Setting $\theta=t_0-t$, we see that $\sigma$ is affine with respect to $t$ if and only if, for all $(\theta,x)\in D_X$,
 $$\bigl\langle d(g-\sigma_{t_0})\circ\Phi_X(\theta,x),X\circ\Phi_X(\theta,x)\bigr\rangle
   =\bigl\langle d(g-\sigma_{t_0})(x),X(x)\bigr\rangle\,,
 $$
which can also be written as
 $${\mathcal L}(X)(g-\sigma_{t_0})\bigl(\Phi_X(\theta,x)\bigr)
   ={\mathcal L}(X)(g-\sigma_{t_0})(x)\,,
 $$
where ${\mathcal L}(X)(g-\sigma_{t_0})$ is the Lie derivative of the function $g-\sigma_{t_0}$ 
with respect to the vector field $X$. This equality expresses the fact that ${\mathcal L}(X)(g-\sigma_{t_0})$
is constant along each integral curve of $X$. 
When the preceding condition is satisfied, the function $\sigma$ can be uniquely extended into a 
function, affine with respect to the variable $t$, 
defined on $\RR\times M$, and one easily check that it satisfies condition $(*)$ of the statement above for all solutions $\varphi:I\to M$ of the differential equation determined by $X$
and all $t\in I$.
\end{proof}

\subsection{Hamiltonian and conformally Hamiltonian vector fields}\label{confhamfields}

\begin{definition}
Let $(M,\omega)$ be a smooth symplectic manifold and $H:M\to\RR$ a smooth function. The unique vector field
$X_H$ such that $i({X}_H)\omega=-dH$
is called the \emph{Hamiltonian vector field} associated to $H$, and $H$ is called a \emph{Hamiltonian} for $X_H$.
Let $g:M\to \RR\backslash\{0\}$ be a nowhere vanishing function. The vector field $Y=g\,{X}_H$ will be called a
\emph{conformally Hamiltonian vector field}, with $H$ as \emph{Hamiltonian} and $g$ as \emph{conformal factor}.
The vector field $Y$ satisfies
 $$i(Y)\omega=gi({X_H})\omega=-g\,dH\,.$$
\end{definition} 

\begin{theorem} \label{result1} 
Let $(M,\omega_1)$ be a symplectic manifold, $H:M\to\RR$ a smooth Hamiltonian, 
$X$ the associated Hamiltonian vector field.
We assume that $X$ is complete; in other words, its flow $\Phi_X$ is defined on the whole of $\RR\times M$.
Let $M^0$ be an open dense subset of $M$, $g:M^0\to\RR\backslash\{0\}$ be a smooth, nowhere vanishing function 
and $Y=gX$ be the conformally Hamiltonian vector field on $M^0$, with Hamiltonian $H$ and conformal factor $g$.
Its flow will be denoted by $\Phi_Y$.
Let $\sigma:\RR\times M^0\to \RR$ be a smooth function such that for each maximal solution $\varphi$ of the
differential equation determined by $Y$,
 $$
 \frac{d\sigma\bigl(t,\varphi(t)\bigr)}{dt}=g\bigl(\varphi(t)\bigr)\,.
 $$
We assume that there exists on $M^0$ another symplectic form $\omega_2$ such that 
 $$i(Y)\omega_2= -dH\,.$$ 
In other words, the vector field $Y$ is both Hamiltonian with respect to $\omega_2$ with $H$ as Hamiltonian and
conformally Hamiltonian with respect to $\omega_1$ with the same $H$ as Hamiltonian and
with $g$ as conformal factor.  

Under these assumptions, the map 
 $$\Xi:M^0\to M\,,\quad x\mapsto\Xi(x)=\Phi_X\bigl(-\sigma(0,x),x\bigr)$$ 
is a symplectic diffeomorphism from $(M^0,\omega_2)$ onto an open subset of $(M,\omega_1)$, equivariant with 
respect to the flow of $Y$ on $M^0$ and the flow of $X$ on $M$, that is
 $$\Xi^*\omega_1=\omega_2
 $$
and for each $(t,x)$ in the open subset of $\RR\times M^0$ on which $\Phi_Y$ is defined
 $$\Phi_X\bigl(t,\Xi(x)\bigr)=\Xi\bigl(\Phi_Y(t,x)\bigr)\,.$$ 
\end{theorem}

\begin{proof}
The maximal integral curve of the differential equation determined by $Y$, which takes the value $ x_0$ for $t=t_0$, is
 $$t\mapsto\Phi_Y(t-t_0, x_0)\,.$$
The same geometric curve in $M^0$, parametrized by $s=\sigma(t,x)$ instead of $t$, is an integral curve of
the differential equation determined by the vector field $X$. The values of the parameter $s$ which
correspond to $(t_0,x_0)$ and to
$\bigl(t,\Phi(t-t_0,x_0))$ are, respectively,
 $$s_0=\sigma(t_0,x_0)\quad\hbox{and}\quad
   s=\sigma\bigl(t, \Phi(t-t_0,x_0)\bigr)\,.$$  
Since $\Phi_X$ is the flow of the vector field $X$, we have
 $$\Phi_Y(t-t_0,x_0)=\Phi_X\Bigl(\sigma\bigl(t,\Phi(t-t_0,x_0)\bigr)-\sigma(t_0,x_0)\,,\,x_0\Bigr)\,.$$ 
Let $\Upsilon:\RR\times M^0\to\RR\times M$ be the map
 $$(t,x)\mapsto \Upsilon(t,x)=\Bigl(t,\Phi_X\bigl(t-\sigma(t,x),x\bigr)\Bigr)\,.$$
We are going to prove that for all $t, t_0, x_0$ such that $\Phi_Y(t-t_0,x_0)$ is defined,
 $$\Upsilon\bigl(t,\Phi_Y(t-t_0,x_0)\bigr)=\Bigl(t, \Phi_X\bigl(t-\sigma(t_0,x_0),x_0\bigr)\Bigr)\,.$$
The above formula expresses the fact that $\Upsilon$ maps the graph of the integral curve $t\mapsto\Phi_Y(t-t_0,x_0)$
of the vector field $Y$ which takes the value $x_0$ for $t=t_0$, into the graph of the integral curve
$t\mapsto\Phi_X\bigl(t-\sigma(t_0,x_0),x_0\bigr)$ which takes the value $x_0$ for $t=\sigma(t_0,x_0)$. 
Observe that $\Upsilon$ associates, to the point of the integral curve of $Y$ reached for the value $t$ of the parameter, the point of the integral curve of $X$ reached for the same value $s=t$ of the parameter.

Replacing $x$ by $\Phi_Y(t-t_0, x_0)$ in the formula which defines $\Upsilon$, we get
$$\Upsilon\bigl(t,\Phi_Y(t-t_0,x_0)\bigr)=\Biggl(t,\Phi_X\Bigl(t-\sigma\bigl(t,\Phi_Y(t-t_0,x_0)\bigr),\Phi_Y(t-t_0,x_0) 
 \Bigr)\Biggr)\,.$$ 
But we have shown that
$$\Phi_Y(t-t_0,x_0)=\Phi_X\Bigl(\sigma\bigl(t,\Phi_Y(t-t_0,x_0)\bigr)-\sigma(t_0,x_0)\,,\,x_0\Bigr)\,.$$
Therefore
 \begin{align*}
 &\Phi_X\Bigl(t-\sigma\bigl(t,\Phi_Y(t-t_0,x_0)\bigr),\Phi_Y(t-t_0,x_0) 
 \Bigr)\\
 &\quad=\Phi_X\Biggl(t-\sigma\bigl(t,\Phi_Y(t-t_0,x_0)\bigr),
 \Phi_X\Bigl(\sigma\bigl(t,\Phi_Y(t-t_0,x_0)\bigr)-\sigma(t_0,x_0)\,,\,x_0\Bigr) 
 \Biggr)\\
 &\quad=\Phi_X\bigl(t-\sigma(t_0,x_0),x_0\bigr)\,.
 \end{align*} 
We have proven that 
$\Upsilon\bigl(t,\Phi(t-t_0,x_0)\bigr)=\Bigl(t, \Psi\bigl(t-\sigma(t_0,x_0),x_0\bigr)\Bigr)$.
Since $Y$ is a Hamiltonian vector field on $(M^0,\omega_2)$, with $H$ as Hamiltonian, the kernel of the closed $2$-form on $\RR\times M^0$
$$\widetilde \omega_2=\omega_2-dH\wedge dt$$
is the rank-one completely integrable distribution generated by the nowhere vanishing vector field
 $$\frac{\partial}{\partial t}+ Y\,,
 $$
where $t$ denotes the coordinate function on the factor $\RR$. We have indeed
 $$i\left(\frac{\partial}{\partial t}+ Y\right)(\omega_2-dH\wedge dt)=dH-dH=0\,.
 $$
Similarly, since $X$ is a Hamiltonian vector field on $(M,\omega_1)$, with $H$ as Hamiltonian, the kernel of the closed 
$2$-form on $\RR\times M$
$$\widetilde \omega_1=\omega_1-dH\wedge ds$$
is the rank-one completely integrable distribution generated by the nowhere vanishing vector field
 $$\frac{\partial}{\partial s}+ X\,,
 $$
where the coordinate function on the factor $\RR$ is now denoted by $s$. We have indeed
 $$i\left(\frac{\partial}{\partial s}+ X\right)(\omega_1-dH\wedge ds)=dH-dH=0\,.
 $$
We recall that $\Upsilon:\RR\times M^0\to\RR\times M$ maps injectively each leaf of the foliation of $\RR\times M^0$
into a leaf of the foliation of $\RR\times M$. On the manifold $\RR\times M^0$, the 2-forms
$\widetilde\omega_2$ and $\Upsilon^*\widetilde\omega_1$ both have the same kernel, since their kernels
determine the same foliation. Each of these $2$-forms is therefore the product of the other by a nowhere vanishing function. This function is in fact the constant $1$, because 
 $$\Upsilon^*s=t\,,\quad, \Upsilon^*H=H\,,\quad\hbox{so}\quad\Upsilon^*(dH\wedge ds)=dH\wedge dt\,,$$
and we have 
$\Upsilon^*\widetilde{\omega_1}=\widetilde{\omega_2}$.

Restricted to $\{0\}\times M^0$, the map $\Upsilon$ becomes
 $$(0,x)\mapsto\Upsilon(0,x)=\Bigl(0, \Psi\bigr(-\sigma(0,x),x\bigr)\Bigr)=\bigl(0,\Xi(x)\bigr)\,.$$
Since $\omega_2$ and $\omega_1$ are the forms induced, respectively, by  $\widetilde\omega_2$ 
on $\{0\}\times M^0$ and by $\widetilde\omega_1$ on $\{0\}\times M$, we have
 $$\Xi^*\omega_1=\omega_2\,.$$ 
The $2$-forms $\omega_1$ and $\omega_2$ being nondegenerate, that proves that the map $\Xi$
is open. But its geometric interpretation proves that $\Xi$ is injective: for each $x_0\in M^0$,
$\Xi(x_0)$ is the point of $M$ reached, for the value $s=0$ of the parameter $s$, by the integral curve $\psi$ of $X$ whose value for $s=s_0=\sigma(0,x_0)$  is $\psi(s_0)=x_0$; it always exists, because $X$ is assumed to be complete. 
Being open and injective, $\Xi$ is a diffeomorphism of $M^0$ onto an open subset of $M$, and more precisely, 
since $\Xi^*\omega_1=\omega_2$, a symplectic diffeomorphism of
$(M^0,\omega_2)$ onto an open subset of $(M,\omega_1)$.

We have seen that $\Upsilon$ is equivariant with respect to the flows of the vector fields
$\displaystyle\frac{\partial}{\partial t}+Y$ on $\RR\times M^0$ and
$\displaystyle\frac{\partial}{\partial s}+X$ on $\RR\times M$. By projection on the second factor, we see that
$\Xi$ is equivariant with respect to the flows of $Y$ on $M^0$ and $X$ on $M$.
\end{proof}

\section{Application to the Kepler problem}\label{kepler}

\subsection{The equations of motion of the Kepler problem}\label{equationsofmotion}
In the physical space, mathematically modelled by an Euclidean affine 3-dimensional space  $\mathcal E$, with associated Euclidean vector space
$\vect{\mathcal E}$, let $P$ be a material point of mass $m$  subjected to the gravitational field created by an
attractive centre $O$. The physical time is mathematically modelled by a real independent variable $t$. Let
$$\vect r = \vect{OP}\,;\quad r=\Vert \vect r\Vert\,;
  \quad \vect p = m\frac{d\vect r}{dt}\,;\quad p=\Vert\vect p\Vert\,.$$
The force $\vect f$ exerted on $P$ is
 $$\vect f = -\frac{km\vect r}{r^3}\,,$$
where $k$ is the constant which characterizes the acceleration field created by $O$.  
\par\smallskip
The equations of motion are
 $$\frac{d\vect r}{dt} =\frac{\vect p}{m}\,,\quad
 \frac{d\vect p}{dt}= -\frac{km\vect r}{r^3}\,.$$
The second equation above becomes singular when $r=0$, so we will assume that the massive point $P$ evolves in
${\mathcal E}\backslash\{O\}$. That space will be called the \emph{configuration space} of the Kepler problem.

\subsection{First integrals}\label{firstintegrals}
The above equations of motion are those determined by a Hamiltonian, time-independent vector field on a symplectic 
6-dimensional manifold, called the \emph{phase space} of the Kepler problem: it is the cotangent bundle
$T^*\bigl({\mathcal E}\backslash\{O\}\bigr)$ to the configuration space. 
We will identify the tangent and cotangent bundles by means of the Euclidean scalar product. The 
Hamiltonian, whose physical meaning is the total energy (kinetic plus potential) of $P$, is
 $$E=E(\vect r, \vect p)=\frac{p^2}{2m}-\frac{mk}{r}\,.$$
The energy $E$ is a first integral of the motion, \emph{i.e.}, it is constant along each 
integral curve. 
\par\smallskip
The group ${\rm SO}(3)$ acts on the configuration space ${\mathcal E}\backslash\{O\}$ by rotations around $O$; 
the canonical lift of that action to the cotangent bundle leaves invariant the Liouville $1$-form and its
exterior derivative, the canonical symplectic form of $T^*\bigl({\mathcal E}\backslash\{O\}\bigr)$. It also leaves
invariant the Hamiltonian $E$. The lifted ${\rm SO}(3)$-action is Hamiltonian, and the corresponding momentum map is an
${\frak so}(3)^*$-valued first integral. Since the Euclidean vector space $\vect{\mathcal E}$ is of dimension $3$,
once an orientation of that space has been chosen, we can identify the Lie algebra ${\frak so}(3)$ and its dual space
${\frak so}(3)^*$ with the vector space $\vect{\mathcal E}$ itself. With that identification the bracket in the Lie algebra
${\frak so}(3)$ becomes the \emph{vector product}, denoted by $(\vect u, \vect v)\mapsto\vect u\times \vect v$; 
the coupling between spaces in duality  becomes the \emph{scalar product} 
$(\vect u, \vect v)\mapsto\vect u\cdot \vect v$; the momentum map of the Hamiltonian
action of ${\rm SO}(3)$ is (up to a sign change) the well known \emph{angular momentum}
 $$\vect L=\vect r\times\vect p\,.$$
In addition to the energy $E$ and the angular momentum $\vect L$, the equations of motion of the Kepler problem have
as a first integral the \emph{eccentricity vector} $\vect{\varepsilon}$
discovered by \emph{Jakob Herman (1678--1753)}~\cite{herman,bernoulli}, improperly called the \emph{Laplace vector}, or the 
\emph{Runge-Lenz vector}, whose origin remained mysterious for a long time:
$$\vect \varepsilon=-\frac{\vect r}{r}+\frac{\vect p\times\vect L}{m^2k}
   = \left(\frac{p^2}{m^2k}-\frac{1}{r}\right)\,\vect r
 -\frac{\vect p\cdot\vect r}{m^2k}\,\vect p\,.
$$
Using the well-known fact (the \emph{first Kepler's law}) that the orbit (the curve described by 
$P$ as a function of time) 
is a conic section with $O$ as one of its foci, the eccentricity vector has a very simple geometric meaning:
it is a dimensionless vector, parallel to the major axis of the orbit, directed from the attracting centre 
towards the perihelion, of length 
numerically equal to the eccentricity of the orbit.

\subsection{The hodograph}\label{hodograph} 
We will not recall in full how the equations of motion of the Kepler
problem can be solved, since it is done in several excellent texts~\cite{anosov, cushmanbates, feynman, milnor}. 
Let us however briefly indicate the proof, due to William Rowan Hamilton (1805--1865)~\cite{hamilton}, 
of an important fact: the hodograph\footnote{The \emph{hodograph} of the motion of a particle moving in an affine space
$\mathcal E$ is the curve, drawn in the associated vector space $\vect{\mathcal E}$, by the velocity vector $\vect v$ of 
the particle, as a function of time.}  
of each motion of the Kepler problem is a circle or an arc of a circle.

Let us first look at solutions for which at the initial time $t_0$, $\vect r$ and $\vect p$ are not 
collinear, so $\vect L(t_0)\neq 0$. Since $\vect L$ is time-independent, $\vect r$
and $\vect p$ are never collinear.
We choose an orthonormal frame, positively oriented, with $O$ as origin and unit vectors 
$\vect{e_x}$, $\vect{e_y}$ and $\vect{e_z}$ with $\vect L$ parallel to 
$\vect{e_z}$. We have
 $$\vect L= L \vect{e_z}\,,\quad\hbox{where $L$ is a constant.}$$ 
The vectors $\vect r$ et $\vect p$ remain for all times parallel to the plane $xOy$. Let $\theta$ 
be the polar angle made by $\vect r$ with $\vect{e_x}$. We have
\begin{equation*}
 \begin{split}
  {\vect r}&=r\cos\theta \vect{e_x}+r \sin\theta \vect{e_y}\,,\\
  {\vect p}&=m\left(\frac{dr}{dt}\cos\theta-r\frac{d\theta}{dt}\sin\theta\right)\vect{e_x}
  +m\left(\frac{dr}{dt}\sin\theta + r\frac{d\theta}{dt}\cos\theta\right)\vect{e_y}\,\\
  \vect L&=mr^2\frac{d\theta}{dt}\vect{e_z}\,.
 \end{split}
\end{equation*}
Therefore   
 $$L=mr^2\frac{d\theta}{dt}=\hbox{Constant}\,.$$
This is the \emph{second Kepler law}, also called \emph{law of areas}, since
$\displaystyle\frac{L}{2m}$ is the area swept by the straight segment $OP$ 
during an unit time, with sign $+$ if $\theta$ increases with time and $-$ if it decreases. Observe that
$\theta$'s variation is strictly monotonic. So we can take $\theta$ as independent variable, instead of time $t$.
We may write
 $$\frac{d\vect p}{d\theta}=\frac{d\vect p}{dt}\,\frac{dt}{d\theta}= 
 -\frac{m^2k}{L}(\cos\theta\vect{e_x}+\sin\theta\vect{e_y})\,. 
 $$
This ordinary differential equation for the unknown $\vect p$, which no more involves $\vect r$, can be readily integrated:
 $$\vect p=\frac{m^2k}{L}(-\sin\theta \vect{e_x}+\cos\theta\vect{e_y})+\vect c\,,$$
where $\vect c$ is a (vector) integrating constant.
We will choose $\vect{e_y}$ such that
 $$\vect c = c \vect{e_y}\,,$$
where $c$ is a numeric constant of the same sign as $L$.

With $O$ as origin let us draw two vectors in the plane $xOy$, the first one (constant) being equal to 
$\vect c$, and the second one (which varies with $\theta$) equal to $\vect p$. 
The end point of that second vector moves on a circle whose centre is the end point of the vector equal to 
$\vect c$, and whose radius is
 $${\mathcal R}=\frac{m^2k}{\vert L\vert}\,.
 $$
This circle (or part of a circle) is (up to multiplication by $m$) the \emph{hodograph} of the
Kepler problem.

A short calculation leads to the following very simple relation between the energy $E$ of a motion, the
radius $\mathcal R$ of its hodograph and the distance $\vert c\vert$ from the attracting centre $O$ to the centre of the 
hodograph: 
 $$2mE=c^2-{\mathcal R}^2\,.$$
Observe that the right-hand side $c^2-{\mathcal R}^2$ is the \emph{power}\footnote{In plane Euclidean geometry,
the \emph{power}\label{power} of a point $O$ with respect to a circle $\mathcal C$ is the real number $\vect{OA}.\vect{OB}$, where $A$
and $B$ are the two intersection points of $\mathcal C$ with a straight line $\mathcal D$ through $O$. That number does not depend
on $\mathcal D$ and is equal to $\Vert \vect{OC}\Vert^2-{\mathcal R}^2$, where $C$ is the centre and $\mathcal R$ 
the radius of $\mathcal C$.} of $O$ with respect to the hodograph.

\subsection{The Levi-Civita parameter}\label{levi-civitaparameter}
Let $\sigma$ be the function, defined on the product with $\RR$ of the phase space of the Kepler problem,
 $$\sigma(t,\vect r, \vect p)=\frac{1}{mk}(\vect p\cdot\vect r-2 E(\vect r, \vect p)t\bigr)\,.$$
A short calculation using the equations of motion shows that for any solution 
$t\mapsto \bigl(\vect{r(t)}, \vect{p(t)}\bigr)$ of the Kepler problem,
 $$\frac{d\sigma\bigl(\vect{r(t)},\vect{p(t)}\bigr)}{dt}=\frac{1}{r(t)}\,.$$ 
With $s(t)=\sigma\bigl(\vect{r(t)},\vect{p(t)})$ as the new independent variable, instead of the time $t$, 
the equations of motion become
  $$\frac{d\vect{r(s)}}{ds}=\frac{r(s)\vect{p(s)}}{m}\,,\quad
    \frac{d\vect{p(s)}}{ds}=-\frac{mk\vect{r(s)}}{r^2(s)}\,.
  $$
I will call $s$ the \emph{Levi-Civita parameter}: it was introduced, as far as I know for the first time, by Tullio~Levi-Civita
in~\cite{levicivita}, by the differential relation which expresses $ds$ as a fuction of $dt$ and $r(t)$. The integrated formula
which gives $\sigma(t,\vect r,\vect p)$ is in the paper~\cite{souriautorino} by Jean-Marie Souriau, and in Exercise 8, chapter II 
of the book by R.~Cushman and L.~Bates~\cite{cushmanbates}. This formula has probably been 
known before for a long time in the Celestial Mechanics community, but I do not know who found 
it for the first time.
\par\smallskip

With the Levi-Civita parameter as new independent variable, the system is no longer Hamiltonian, but rather \emph{conformally Hamiltonian}. 
The Levi-Civita parameter will be used in subsection~\ref{symplecticdiff} in a slightly different context: we will define a diffeomorphism $S$ from the phase space of the Kepler problem restricted to negative 
(resp. positive) values of the energy, onto an open dense subset of the cotangent bundle to a three-dimensional sphere
(resp., to one sheet of a two-sheeted three-dimensional hyperboloid). On this new phase space equipped with its canonical symplectic form, the image of the vector field which determines the equations of motion of the Kepler problem will be conformally Hamiltonian instead of Hamiltonian, while the image of the vector field transformed by the use of the Levi-Civita parameter as independent variable will be Hamiltonian. 
 
\subsection{The Gy\"orgyi-Ligon-Schaaf symplectic diffeomorphism} \label{glsdiffeo}
G.~Gy\"orgyi~\cite{gyorgyi} gave the expression of a symplectic diffeomorphism from the phase space of 
the Kepler problem restricted to the negative (resp. positive) values of the energy, onto an open subset
of the cotangent bundle to a three-dimensional sphere (resp. two-sheeted hyperboloid). He obtained
this diffeomorphism in two steps. First, following the ideas of Fock~\cite{fock}, 
employed two years later by Moser in his well known paper~\cite{moser}, he arranged the components of  the first integrals
$E$, $\vect L$ and $\vect\varepsilon$ of the Kepler problem into a $4\times4$ matrix, and used that matrix to define two
$4$-dimensional vectors $\vect\rho$ and $\vect\pi$ (formulae 2.15 and 2.16 of Gy\"orgyi's 1968 paper~\cite{gyorgyi}).
I believe that the map $(\vect r, \vect p)\mapsto (\vect\rho,\vect\pi)$ obtained by Gy\"orgyi after difficult to follow calculations, is 
the map, maybe rescaled, I have called $S^{-1}$ in Subsection~\ref{smap} below. Second, in Section entiteled \lq\lq Getting back the time $t$; 
Bacry's generators\rq\rq\ of his 1968 paper, Gy\"orgyi composes the map $(\vect r, \vect p)\mapsto (\vect\rho,\vect\pi)$ with a suitably chosen map
deduced from the flow of the Kepler vector field. He did not prove that the map he finally obtained is symplectic. This result was proven
by Ligon and Schaaf~\cite{ligonschaaf} who, eight years later, rediscovered the same diffeomorphism. Twenty years later, 
R.~H.~Cushman and J.~J.~Duistermaat~\cite{cushmanduistermaat} discussed the properties of that diffeomorphism, 
gave new proofs of its remarkable properties and stated the nice characterization indicated in the Introduction.
\par\smallskip

For negative values of the energy, the Gy\"orgyi-Ligon-Schaaf diffeomorphism is the map which, 
to each element $(\vect r, \vect p)$ of the phase space of the Kepler problem such that $E(\vect r, \vect p)<0$, 
associates the element $\bigl((\xi_0,\vect \xi),(\eta_0,\vect\eta)\bigr)$ of $\RR^4\times\RR^4$
\begin{align*}
 \xi_0&=\frac{\sqrt{-2mE(\vect r, \vect p)}}{mk}\,\vect r.\vect p\sin\varphi
       +\left(\frac{rp^2}{mk}-1\right)\cos\varphi\,,\\
 \vect \xi&=\left(\frac{\vect r}{r}-\frac{\vect r.\vect p}{mk}\,\vect p\right)\sin\varphi
       +\frac{\sqrt{-2mE(\vect r, \vect p)}}{mk}\,r\vect p\cos\varphi\,,\\
 \eta_0&=-\vect r.\vect p\cos\varphi 
       +\frac{mk}{\sqrt{-2mE(\vect r,\vect p)}}\left(\frac{rp^2}{mk}-1 \right)\sin\varphi\,,\\       
 \vect\eta&=-\frac{mk}{\sqrt{-2mE(\vect r, \vect p)}}\left(\frac{\vect r}{r}-\frac{\vect r.\vect p}{mk}\,
       \vect p\right) \cos\varphi +r\vect p\sin\varphi\,.\\
\end{align*}
In these formulae $\varphi$ is the angle given, as a function of $(\vect r, \vect p)$, by       
 $$\varphi(\vect r,\vect p)=\frac{\sqrt{-2mE(\vect r,\vect p)}}{mk\vect r.\vect p}\,.$$
The quantities $\vect \xi=(\xi_1,\xi_2, \xi_3)$ and $\vect\eta=(\eta_1,\eta_2, \eta_3)$ are vectors of
$\RR^3$, while $(\xi_0,\vect\xi)=(\xi_0,\xi_1,\xi_2, \xi_3)$ and $(\eta_0,\vect\eta)=(\eta_0,\eta_1,\eta_2,\eta_3)$
are vectors of $\RR^4$, such that
\begin{align*}
 \xi_0^2+\Vert\vect\xi\Vert^2&=\xi_0^2+\xi_1^2+\xi_2^2+\xi_3^2=1\,,\\
 \xi_0\eta_0+\vect{\mathstrut\xi}.\vect{\mathstrut\eta}&=\xi_0\eta_0+\xi_1\eta_1+\xi_2\eta_2+\xi_3\eta_3=0\,,\\
 \eta_0^2+\Vert\vect\eta\Vert^2&=\eta_0^2+\eta_1^2+\eta_2^2+\eta_3^2>0\,.
\end{align*}
In other words, $(\xi_0,\vect\xi)$ is a point of the sphere $S^3$ of radius 1 centered at the origin, 
embedded in $\RR^4$ (endowed with its usual Euclidean metric), and $(\eta_0,\vect\eta)$ is a non-zero vector tangent to that sphere at point $(\xi_0,\vect\xi)$. Using the usual scalar product of $\RR^4$ to identify vectors and
covectors, we may consider $\bigl((\xi_0,\vect\xi),(\eta_0,\vect\eta)\bigr)$ as an element of the cotangent bundle to the sphere $S^3$ minus its zero section.

For positive values of the energy, the Gy\"orgyi-Ligon-Schaaf diffeomorphism is given by similar formulae, the trigonometric functions $\sin$ and $\cos$ being replaced by the hyperbolic functions $\sinh$ and $\cosh$ and the
Euclidean metric on $\RR^4$ by the Lorentz pseudo-Euclidean metric.

We will show in what follows that the Gy\"orgyi-Ligon-Schaaf symplectic diffeomorphism can be easily obtained, and
that all its properties can be proven, by an adaptation of Moser's method for the regularization of the Kepler problem and application of Theorem~\ref{result1}. 
That method rests on the fact that hodographs of the Kepler problem are circles or arcs of circles, and that
the stereographic projection maps circles into circles. We do not need to know in advance that the eccentricity
vector $\vect\varepsilon$ is a first integral of the Kepler problem, since this property is an easy consequence of Moser's method.

\subsection{Stereographic projection}\label{stereoproj}
In 1935, V.~A.~Fock~\cite{fock} used an inverse stereographic projection for the determination of the energy levels
of hydrogen atoms in quantum mechanics. More recently J.~Moser~\cite{moser} used the same idea for the regularization of the Kepler problem. His method yields a correspondence between oriented hodographs of motions of the Kepler problem
with a fixed, negative energy, and great oriented circles of a three-dimensional sphere. 
As explained by D.~V.~Anosov~\cite{anosov} and J.~Milnor~\cite{milnor}, the method can easily be adapted for motions 
with a fixed, positive energy: instead of a three-dimensional sphere, one should use one sheet of a two-sheeted 
three-dimensional hyperboloid of revolution, and one obtains a correspondence between oriented hodographs 
of motions of the Kepler problem
with the chosen positive energy, and connected components of great hyperbolas drawn on that sheet of hyperboloid
(that means intersections of that sheet of hyperboloid with a plane containing the symmetry centre).

To deal with cases $E<0$ and $E>0$ simultaneously, we introduce the auxiliary quantity
 \begin{equation*}
  \zeta =
  \begin{cases}
  1 &\hbox{if}\quad E<0\,,\cr
  -1 &\hbox{if}\quad E>0\,.\cr
  \end{cases} 
 \end{equation*}
Let $(O, \vect{e_x},\vect{e_y},\vect{e_z})$ be an orthonormal frame of $\mathcal E$.  with $O$ as origin.
We add to the basis $(\vect{e_x},\vect{e_y},\vect{e_z})$ of $\vect{\mathcal E}$ a fourth vector $\vect{e_h}$ 
and we denote by  $h$ the corresponding coordinate. So we get a 4-dimensional affine space 
$\mathcal F$. The physical space $\mathcal E$ will be identified with the affine subspace of 
$\mathcal F$ determined by the equation $h=0$.

In $\mathcal F$, for each $\rho>0$, let $Q_\rho$ be the quadric defined by
 $$h^2+\zeta(x^2+y^2+z^2)=\rho^2\,.$$
It is a sphere if $\zeta=1$, a two-sheeted hyperboloid if $\zeta=-1$.  
We will see that there is a value of $\rho$ particularly suited for each value of the energy $E$.

\begin{figure}[htp]
\begin{center}

  \includegraphics{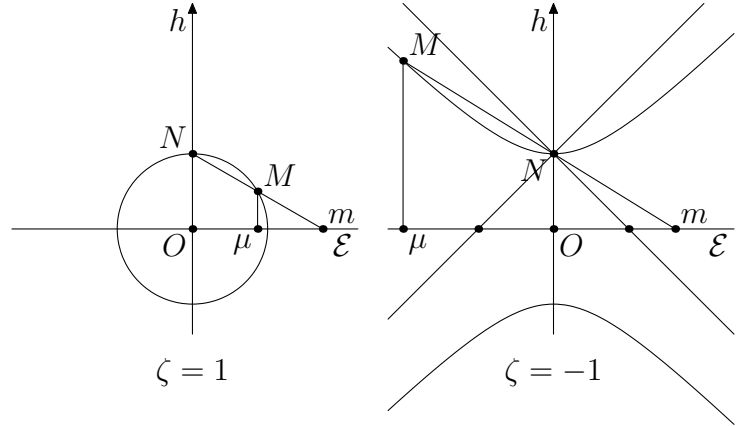}\\
  \caption{Stereographic projection and its generalization}\label{projstereo}
  \end{center}
\end{figure}

Let $N$ be the point with coordinates $(x=y=z=0,\quad h=\rho)$. The  \emph{stereographic projection}
(usual if $\zeta=1$, generalized if $\zeta=-1$) of the quadric $Q_\rho$ minus point $N$ on the space $\mathcal E$
is the map which associates, with each point $M\in Q_\rho\backslash\{N\}$, the intersection point $m$  
of the straight line which joins $N$ and $M$, with  $\mathcal E$. See Figure~\ref{projstereo}.
If $\zeta=1$  that map is a diffeomorphism from $Q_\rho\backslash\{N\}$ onto $\mathcal E$.
If $\zeta=-1$ it is a diffeomorphism from $Q_\rho\backslash \{N\}$ onto the open subset of $\mathcal E$ complementary to the 2-sphere of centre $O$ and radius $\rho$. The upper sheet ($h>0$) of the hyperboloid (minus point $N$) is mapped onto the outside of that sphere and the lower sheet ($h<0$) onto its inside.

\subsection{Cotangent lift}\label{cotangentlift}
Since the (maybe generalized) stereographic projection is a diffeomorphism, it can be uniquely lifted to the cotangent
bundles in such a way that the pull-back of the Liouville form on $T^*{\mathcal E}$ is equal to the Liouville form on 
$T^*(Q_\rho\backslash\{N\})$.
We apply this construction with $\vect{Om}=\vect p$, the linear momentum of a Keplerian motion of energy $E\neq 0$.

Each $1$-form on $\mathcal E$ can be written $\vect r\cdot d\vect{Om}$, where $\vect r$ is a vector field on
 $\mathcal E$, or, since $\vect{Om}=\vect p$, 
 $$\vect r\cdot d\vect p=r_x\,dp_x+r_y\,dp_y+r_z\,dp_z\,.$$
Recall that $\zeta=1$ if $E<0$ and $\zeta =-1$ if $E>0$. The equation of $Q_\rho$ is
 $$\vect{OM}\cdot \vect{OM}= x^2+y^2+z^2+\zeta h^2=\zeta \rho^2\,.$$
Each $1$-form on $Q_\rho$ can be written
 $$\vect W\cdot d\vect{OM}= W_x\, dx+W_y\,dy+W_z\,dz+\zeta W_h\,dh\,,$$ 
where $\vect W$ is a vector field defined on $Q_\rho$, tangent to that quadric. Therefore, 
 $$\vect W\cdot \vect{OM}=W_x x +W_y y +W_z z+\zeta W_h h=0\,.$$

The cotangent lift of the (maybe generalized) stereographic projection will be denoted by $S_\rho$. It maps 
each pair $(M, \vect W)$ made by a point $M\in Q_\rho\backslash\{N\}$ and a vector
$\vect W$ tangent to $Q_\rho$ at that point, onto a pair $(\vect r, \vect p)$ of two vectors of $\vect{\mathcal E}$, in
such a way that
 $$\vect W\cdot d\vect{OM}=\vect r\cdot d\vect p\,.$$
The formulae which give $(\vect p, \vect r)$ as functions of $(\vect{OM},\vect W)$ 
are easily obtained after some calculations:
\begin{equation*}
\left\{
 \begin{aligned}
  \vect p&=\frac{\rho}{\rho-h}\,\vect {O\mu}\,,\quad\hbox{with}\quad h^2+\zeta\Vert \vect{O\mu}\Vert^2=\rho^2\,,\\
  \vect r&=\frac{\rho-h}{\rho}\,\vect {W_3}+\frac{W_h}{\rho}\,\vect{O\mu}\,,
 \end{aligned}
\right.
\end{equation*}
where we have set
 $$\vect{OM} =\vect{O\mu}+h\vect{e_h}\,,\quad \vect W=\vect{W_3}+W_h\vect{e_h}\,,$$ 
 The formulae for the inverse transform are:
\begin{equation*}
\left\{
 \begin{aligned}
  \vect{O\mu}&=\frac{2\rho^2}{\rho^2+\zeta p^2}\vect p\,,\\
  h&=\rho\,\frac{p^2-\zeta\rho^2}{p^2+\zeta\rho^2}\,,
     \quad\hbox{with}\quad p=\Vert\vect p\Vert=\sqrt{p_x^2+p_y^2+p_z^2}\,,\\
  \vect{W_3}&=\frac{\rho^2+\zeta p^2}{2\rho^2}\vect r-\frac{\zeta\vect r.\vect p}{\rho^2}\vect p\,,\\
  W_h&=\frac{\zeta\vect r.\vect p}{\rho}\,.
 \end{aligned}
\right.
\end{equation*}
These formulae prove that if $\vect{p_1}=\vect{Om_1}$ and $\vect{p_2}=\vect{Om_2}$ are two 
collinear vectors in $\vect{\mathcal E}$ such that 
 $$\vect{Om_1}.\vect{Om_2}=-\zeta \rho^2\,,$$
their images, by the inverse stereographic projection, are two points $M_1$ and $M_2$ 
symmetric of each other with respect to $O$. That property may also be proven without calculation
as a consequence of properties of the transformation by inversion. 

We recall that if a straight line through $O$ meets the hodograph of a Keplerian motion at two points
$m_1$ and $m_2$, $\vect{Om_1}.\vect{Om_2}$ is the power (see the second footnote in Subsection~\ref{hodograph}) of $O$ with respect to the hodograph, and is equal to
$2mE$, where $E$ is the energy of the motion. Therefore the above results shows that the inverse streographic projection
maps the hodographs of motions of energy $\displaystyle -\frac{\zeta\rho^2}{2m}$ on curves, drawn 
on the quadric $Q_\rho$, which are symmetric with respect to $O$. These curves are 
great circles of the sphere $Q_\rho$ if $\zeta=1$, and
great hyperbolas of the hyperboloid $Q_\rho$ if $\zeta=-1$. The name \lq\lq great hyperbola\rq\rq\ is used
by analogy with \lq\lq great circle\rq\rq: it means the intersections of $Q_\rho$ with a plane containing 
the symmetry centre $O$. As a consequence, the cotangent lift $S_\rho^{-1}$ of the inverse stereographic projection
maps the hodographs of motions of the Kepler problem of energy
$\displaystyle E=-\frac{\zeta\rho^2}{2m}$ on geodesics of the quadric $Q_\rho$. The hodographs
of motions with a different value of the energy are mapped on circles (if $\zeta=1$) or hyperbolas
(if $\zeta=-1$)  drawn on the quadric $Q\rho$, but which are not \emph{great} circles, or \emph{great} hyperbolas.

So we see that for each value of the energy $E$, there is a particularly well suited value of $\rho$: it is the unique
$\rho>0$ such that
$\displaystyle E=-\frac{\zeta\rho^2}{2m}$, with $\zeta=1$ if $E<0$ and $\zeta=-1$ if $E>0$.

\begin{remark}\label{antisymp}
We have defined the mapping $S_\rho:T^*(Q_\rho\backslash\{N\})\to T^*{\mathcal E}$  so that
 $$(S_\rho)^*(\vect r\cdot d\vect{p})=\vect W\cdot d\vect{OM}\,.$$
Since the symplectic form of  $T^*{\mathcal E}$ is $d(\vect p\cdot d\vect{r})$ and that of $T^*(Q_\rho\backslash\{N\})$ 
$d(\vect W\cdot d\vect{OM})$, the map $S_\rho$ is an \emph{anti-symplectic} diffeomorphism, rather than a symplectic diffeomorphism.
\end{remark} 

The transformed Hamiltonian $(S_\rho)^* E=E\circ S_\rho$ is
$$E\circ S_\rho =\frac{-\zeta\rho^2}{2m} +\frac{\zeta\rho^3}{m(\rho-h)\Vert \vect W\Vert}
 \left(\Vert\vect W\Vert-\frac{km^2}{\rho^2}\right)\,,$$ 
where we have set
 $$\Vert\vect W\Vert=\sqrt{W_x^2+W_y^2+W_z^2+\zeta W_h^2}\,.$$
The expression of $E\circ S_\rho$ shows that $\displaystyle E\left(\vect r,\vect p\right)=\frac{-\zeta\rho^2}{2m}$
if and only if $S_\rho^{-1}\left(\vect r, \vect p\right)=\left(\vect{OM}, \vect W\right)$ is such that
$\displaystyle\Vert \vect W\Vert=\frac{km^2}{\rho^2}$. 
By differentiation we obtain
 $$d(E\circ S_\rho) =\frac{\zeta\rho^3}{m(\rho-h)\Vert \vect W\Vert}\,d
 \left(\Vert\vect W\Vert\right)+\left(\Vert\vect W\Vert-\frac{km^2}{\rho^2}\right)
 d\left(\frac{\zeta\rho^3}{m(\rho-h)\Vert \vect W\Vert}\right)\,.
 $$  
For a motion of energy $\displaystyle E=-\frac{\zeta\rho^2}{2m}$, the second term of the right hand side
of the above equality vanishes, since
$\displaystyle  \left(\Vert\vect W\Vert-\frac{km^2}{\rho^2}\right)=0$, and we see that
$S_\rho^{-1}$ transforms the Hamiltonian vector field of the Kepler problem into a conformally
Hamiltonian vector field on $T^*(Q_\rho\backslash\{N\})$ minus the zero section, with Hamiltonian 
$\Vert \vect W\Vert$ and conformal factor $\displaystyle \frac{-\zeta\rho^3}{m(\rho-h)\Vert \vect W\Vert}$. 
We introduced a minus sign in the conformal factor to account for the fact that $S_\rho^{-1}$ is anti-symplectic
(see Remark~\ref{antisymp}).

\subsection{Motions of zero energy}\label{zeroenergy}
For Keplerian motions of energy zero the stereographic projection should be replaced by
an inversion with pole $O$ and ratio $ l$. Formulae for its cotangent lift are 
\begin{equation*}
\left\{
 \begin{aligned}
  \vect{p}&=\frac{ l}{\Vert\vect{OM}\Vert^2}\,\vect{OM}\,,\\
  \vect r&=\frac{\Vert\vect{OM}\Vert^2}{ l}\,\vect W-\frac{2\vect W\cdot \vect{OM}}{ l}\,\vect{OM}\,.
 \end{aligned}
\right.
\end{equation*}
The correspondence $(\vect p, \vect r)\mapsto(\vect{OM}, \vect W)$ being involutive, the formulae
for the inverse transformation are
\begin{equation*}
\left\{
 \begin{aligned}
  \vect{OM}&=\frac{ l}{p^2}\,\vect p\,,\\
  \vect W&=\frac{p^2}{ l}\,\vect r-\frac{2\vect r\cdot \vect p}{ l}\,\vect p\,.
 \end{aligned}
\right.
\end{equation*}
Denoting by $S_0$ the cotangent lift of the inversion of pole $O$ and ratio $ l$, we easily obtain
the transformed Hamiltonian 
 $$E\circ {S_0}=\frac{ l^2}{2m\Vert \vect W\Vert\,\Vert\vect{OM}\Vert^2}\,
 \left(\Vert \vect W\Vert - \frac{2m^2k}{ l} \right)\,.$$
Now $\vect W$ is a vector of $\vect{\mathcal E}$, with 3 components, so we have set 
 $$\Vert\vect W\Vert=\sqrt{W_x^2+W_y^2+W_z^2}\,.$$
As above this result proves that  $E\left(\vect r,\vect p\right)=0$
if and only if $S_0^{-1}\left(\vect r, \vect p\right)=\left(\vect{OM}, \vect W\right)$ is such that
$\displaystyle\Vert\vect W\Vert=\frac{2km^2}{ l}$. 

By differentiation we obtain
 $$d(E\circ {S_0})=\frac{ l^2}{2m\Vert \vect W\Vert\,\Vert\vect{OM}\Vert^2}\,
 d\left(\Vert \vect W\Vert\right)+\left(\Vert \vect W\Vert - \frac{2m^2k}{ l} \right)
 d\left(\frac{ l^2}{2m\Vert \vect W\Vert\,\Vert\vect{OM}\Vert^2}
 \right)\,.
 $$
For a motion of energy $\displaystyle E=0$, the second term of the right hand side
of the above equality vanishes, since
$\displaystyle  \left(\Vert\vect W\Vert-\frac{2km^2}{ l}\right)=0$, and we see that
$S_0^{-1}$ transforms the Hamiltonian vector field of the Kepler problem into a conformally
Hamiltonian vector field on $T^*({\mathcal E}\backslash\{0\})$ minus the zero section, with Hamiltonian 
$\Vert \vect W\Vert$ and conformal factor $\displaystyle \frac{- l^2}{2m\Vert \vect W\Vert\,\Vert\vect{OM}\Vert^2}$. 
We introduced a minus sign in the conformal factor to account for the fact that $S_0$ is anti-symplectic
(Remark~\ref{antisymp} applies to $S_0$ as well as to $S_\rho$).

\subsection{Infinitesimal symmetries}\label{infinitesimalsymmetries} Let us first recall some properties 
about infinitesimal symmetries. Consider a vector field $X$ on a manifold $M$. Another vector field $Z$
on $M$ is said to be an \emph{infinitesimal symmetry} of the differential equation determined by $X$ 
(or, in short, of $X$) if $[Z,X]=0$; it is said to be a \emph{weak infinitesimal symmetry} of $X$ 
if at each point $x\in M$, $[Z,X](x)$ and $X(x)$ are collinear. If $Z$ is an infinitesimal symmetry or a weak infinitesimal symmetry of $X$
and $h$ a smooth function, $Z$ is a weak infinitesimal symmetry of $hX$, since
 $$[Z, hX]= h[Z,X]+\bigl({\mathcal L}(Z)h\bigr)\,X\,,$$
whichs shows that at each point  $x\in M$, $[Z,hX](x)$ and $X(x)$ are collinear. 
When $X=X_H$ is a Hamiltonian vector field on a symplectic manifold 
$(M,\omega)$, with a function $H$ as Hamiltonian, a vector field $Z$ on $M$ such that ${\mathcal L}(Z)\omega=0$ and 
${\mathcal L}(Z)H=0$ is an infinitesimal symmetry of $X_H$ and a weak infinitesimal symmetry of
$hX_H$ for any smooth function $h$. When in addition $Z$ is Hamiltonian, its 
Hamiltonian is a first integral of $X_H$ and of $gX_H$. Conversely, if $f$ and $g$ are two smooth first integrals 
of the Hamiltonian vector field $X_H$, the vector field $gX_f$, $X_f$ being the Hamiltonian vector field with Hamiltonian $f$,
is an infinitesimal symmetry of $X_H$ and a weak infinitesimal symmetry of  $hX_H$ for any smooth function $h$. We have indeed
 \begin{align*}
 [gX_f,hX_H]&=gh[X_f,X_H]+g\bigl(i(X_f)h\bigr)\,X_H-h\bigl(i(X_H)g\bigr)\,X_f\\
            &=ghX_{\{f,H\}}+g\{f,h\}\,X_H-h\{H,g\}\,X_f=g\{f,h\}\,X_H\,,
 \end{align*}
since $\{f,H\}=\{H,g\}=0$.
\par\smallskip

For the vector field which determines the equations of motion of the Kepler problem, 
we already know a three-dimensional vector space of infinitesimal symmetries: it 
is made by the canonical lifts to the cotangent bundle 
of infinitesimal rotations of the configuration space ${\mathcal E}\backslash\{O\}$ 
around the attractive centre $O$. These canonical lifts are Hamiltonian vector fields, 
and their Hamiltonians (linear combinations of the components of the angular momentum $L$)
are first integrals of the equations of motion. The results obtained in 
Subsections~\ref{cotangentlift} and~\ref{zeroenergy} 
show that other weak infinitesimal symmetries exist. In these subsections, 
for each value $e$ of the energy $E$, we have 
built a symplectic diffeomorphism\footnote{Or rather an anti-symplectic diffeomerphism, 
but by a change of sign it can be transformed into a symplectic diffeomorphism} 
of the phase space of the Kepler problem onto an open subset of another symplectic manifold: 
the cotangent space to a three-dimensional 
sphere when $e<0$, the cotangent space to $\mathcal E$ when $e=0$ and the cotangent space to a three-dimensional
hyperboloid when $e>0$. That sphere, affine space $\mathcal E$ or  hyperboloid
will be called the \emph{new configuration space}, and its cotangent bundle (with the zero section removed) 
the \emph{new phase space}. 
\par\smallskip

For each energy level $e$, the new phase space has a 6-dimensional Lie group of global symmetries: 
${\rm SO}(4)$ when $e<0$, the group ${\rm SE}(3)$ of orientation-preserving isometries (rotations and translations) 
of $\mathcal E$ when $e=0$
and the Lorentz group ${\rm SO}(3,1)$ when $e>0$. The direct image of the Hamiltonian vector field 
of the Kepler problem, restricted to each energy level of the  phase space, 
is equal to the restriction, to the corresponding submanifold of the new phase space, 
of a conformally Hamiltonian vector field, whose Hamiltonian is invariant by the action of the 
$6$-dimensional Lie group of symmetries of the new phase space. Therefore, on each energy level of the phase space 
of the Kepler problem, six linearly independent weak infinitesimal symmetries exist, of which three are the already known 
linearly independent infinitesimal symmetries associated with the cotangent lifts of infinitesimal rotations
of $\mathcal E$ around three non-coplanar axes through $O$. The other three correspond to the additional
symmetries of the considered energy level of the new phase space.  
We will see that they depend smoothly 
on the value of the energy. Therefore, on the phase space of the Kepler problem considered as a whole, 
there exist $6$ linearly independent vector fields, three of them being infinitesimal symmetries and the 
other three weak infinitesimal symmetries. We will see that these six vector fields, which
of course are tangent to each energy level, are Hamiltonian and are all
infinitesimal symmetries of the Kepler vector field. So their Hamiltonians are first integrals of the equations of motion.
\par\smallskip

In Subsection~\ref{energymomentum} we will see that these infinitesimal symmetries 
do not span a Lie algebra, but rather a fibered space in Lie 
algebras~\cite{douadylazard} over the base $\RR$, the value of the energy running over that base, in other words a 
\emph{Lie algebroid}~\cite{cannasweinstein, mackenzie} of a special kind, with a zero anchor map. 
\par\smallskip

Let us recall~\cite{libermannmarle, ortegaratiu} that if $\Phi:G\times M\to M$ is a left action of 
a Lie group $G$ on a manifold~$M$, its canonical lift to the cotangent bundle is a Hamiltonian action
$\widehat\Phi:G\times T^*M\to T^*M$ with an equivariant momentum map
$J:T^*M\to{\mathcal G}^*$ given by 
 $$\langle J(\xi), X\rangle = \bigl\langle \xi, X_M\bigl(\pi_M(\xi)\bigr)\bigr\rangle\,.
 $$
In this formula, $\pi_M:T^*M\to M$ is the canonical projection; $X\in{\mathcal G}$, the Lie algebra of $G$; 
$\xi\in{\mathcal G}^*$, the dual space of $\mathcal G$; and $X_M$ is the fundamental vector field on $M$ 
associated to $X$, defined by
 $$X_M(x)=\frac{d\Bigl(\Phi\bigl(\exp(tX),x\bigr)\Bigr)}{dt}\Biggm|_{t=0}\quad\text{for each}\quad x\in M\,.$$
The momentum map $J$ is a first integral of any Hamiltonian vector field $X_H$ or conformally Hamiltonian
vector field $gX_H$ on $T^*M$ whose Hamiltonian $H$ is invariant under $\widehat \Phi$. 
For each $X\in {\mathcal G}$ the Hamiltonian vector field on $T^*M$
whose Hamiltonian is $\xi\mapsto\langle J(\xi),X\rangle$ is an infinitesimal symmetry of $X_H$.
\par\smallskip
 
We apply these results, with for $\widehat\Phi$ the action on the new phase space of its symmetry group.
As a vector space, the Lie algebra of the symmetry group 
is canonically isomorphic to $\vect{\mathcal E}\times\vect{\mathcal E}$, but with a Lie algebra bracket 
$[\ ,\ ]_e$ which depends on the energy level $e$. That bracket will be determined in Subsection~\ref{energymomentum}. 
Two kinds of infinitesimal symmetries of the new phase space exist. The first kind is made by the fundamental
vector fields which correspond to elements $(\vect{u_1},0)$ of $\vect{\mathcal E}\times\vect{\mathcal E}$. 
These vector fields are the canonical lifts to the cotangent
bundle of the infinitesimal rotation of $\mathcal E$ around the axis through $O$ parallel to $\vect{u_1}$.
The second kind is made by the fundamental vector fields which correspond to elements $(0,\vect{u_2})$ of 
$\vect{\mathcal E}\times\vect{\mathcal E}$. These vector fields are the canonical lifts 
to the cotangent bundle of the following infinitesimal transformations of the new configuration space:

\begin{itemize}

\item{} when $e<0$, the restriction to the sphere $Q_\rho$, with $\rho^2=-2me$, of the infinitesimal rotation of 
$\vect{\mathcal F}=\vect{\mathcal E}\times(\RR\times\vect{e_h})$ in which the plane
spanned by $\vect{u_2}$ and $\vect{e_h}$ rotates, while the vector subspace of $\vect{\mathcal E}$ orthogonal to  
$\vect{u_2}$ remains fixed;    

\item{} when $e=0$, the infinitesimal translation of $\mathcal E$ parallel to $\vect{u_2}$;

\item{} when $e>0$, the restriction to the hyperboloid $Q_\rho$, with $\rho^2=2me$, of the infinitesimal Lorentz transformation 
of $\vect{\mathcal F}=\vect{\mathcal E}\times(\RR\times\vect{e_h})$ in which the plane
spanned by $\vect{u_2}$ and $\vect{e_h}$ is transformed, while the vector subspace of $\vect{\mathcal E}$ orthogonal to  
$\vect{u_2}$ remains fixed.    
\end{itemize} 
The space $\vect{\mathcal E}\times\vect{\mathcal E}$ will be identified 
with its dual space by using the scalar product on each factor.
The momentum map, which will be denoted by $K_e$,
can be written $K_e=(K_{1e},K_{2e})$, each component taking its values in $\vect{\mathcal E}$.
\par\smallskip

In the three cases $e<0$, $e=0$ and $e>0$, the Hamiltonian $H$ has the same expression
 $$(\vect{OM},\vect W)\mapsto H(\vect{OM},\vect W)=\Vert \vect W\Vert\,.$$
It is equivalent to the classical  Hamiltonian 
$\displaystyle \frac{\Vert\vect W\Vert^2}{2}$ of a particle moving freely on $Q_\rho$ 
when $e\neq 0$ or on $\mathcal E$ when $e=0$. It is indeed invariant under the action $\widehat \Phi$.

The first component $K_{1e}$ of the momentum map is
\begin{equation*}
 K_{1e}(\vect{OM},\vect{W})=
 \begin{cases}
  \vect{O\mu}\times\vect{W_3}&\text{when\quad $e\neq 0$,}\\
  \vect{OM}\times\vect{W}&\text{when\quad $e=0$,}
 \end{cases}
\end{equation*} 
and its second component is
\begin{equation*}
 K_{2e}(\vect{OM},\vect{W})=
 \begin{cases}
  \zeta(h\vect{W_3} - W_h\vect{O\mu})&\text{when\quad $e\neq 0$,}\\
  \vect W&\text{when\quad $e=0$.}
 \end{cases}
\end{equation*}
As above, we have set when $e\neq 0$
 $$\vect{OM}=\vect{O\mu}+h\vect{e_h}\,,\quad \vect{W}=\vect{W_3}+W_h\vect{e_h}\,.$$
Let us now come back to the phase space of the Kepler problem. The composed map 
$J_e=K_e\circ S_\rho^{-1}$, which has two components $J_{1e}=K_{1e}\circ S_\rho^{-1}$ and
$J_{2e}=K_{2e}\circ S_\rho^{-1}$, is the momentum map of the Lie algebra action 
of infinitesimal symmetries. Using the formulae which give
$\vect {OM}$ and $\vect W$ as functions of $\vect r$ and $\vect p$, we obtain for the first component, 
in the three cases $e<0$, $e=0$ and $e>0$
\begin{equation*}
 J_{1e}=\vect p\times \vect r\,.
\end{equation*}
Its expression does not depend on $e$. Up to a change of sign, it is the angular momentum $\vect L$. For the second component we obtain
\begin{equation*}
 J_{2e}=
\begin{cases}
 \displaystyle 
 \frac{p^2-\zeta\rho^2}{2\rho}\,\vect r-\frac{\vect r.\vect p}{\rho}\,\vect p&\text{when\quad $e\neq0$,}\\
 \displaystyle
 \frac{p^2}{l}\,\vect r-\frac{2\vect r.\vect p}{l}\,\vect p&\text{when\quad $e=0$.}
\end{cases}
\end{equation*}
Using the value of the energy $\displaystyle E=\frac{p^2}{2m}-\frac{mk}{r}$
and, when $e\neq 0$, the equality $\zeta \rho^2=-2mE$, these formulae become
\begin{equation*}
 J_{2e}=
\begin{cases}
 \displaystyle 
 \frac{m^2k}{\sqrt{-2\zeta m e}}\,\vect\varepsilon&\text{when\quad $e\neq0$,}\\
 \displaystyle
 \frac{2m^2k}{l}\,\vect \varepsilon&\text{when\quad $e=0$,}
\end{cases}
\end{equation*}
where $\vect \varepsilon$ is the eccentricity vector. Since the energy $E(\vect r, \vect p)$ is a first integral, 
we see that by choosing $l=2m^2k$
and, for $e\neq 0$, by multiplying the weak infinitesimal symmetries by the smooth function
$\displaystyle \frac{\sqrt{-2\zeta m E(\vect r, \vect p)}}{m^2k}$, we can arrange things so that
the second component of the momentum map becomes
 $$J_{2e}=\vect\varepsilon\,,
 $$
which no longer depends on the energy level $e$. Therefore, although we 
have defined them separately on each energy level, the three additional weak infinitesimal 
symmetries smoothly depend on the energy level, since they are the Hamiltonian 
vector fields whose Hamiltonians are the components of
$\vect \varepsilon$; they are true (not only weak) infinitesimas symmetries since they are Hamiltonian vector fields whose
Hamiltonians are first integrals of the Kepler equations of motion.

\subsection{The energy-momentum space and map}\label{energymomentum}
The Poisson brackets of the components of the angular momentum $\vect L$ and of the eccentricity vector
$\vect\varepsilon$ in an orthonormal frame $(\vect{e_x},\vect{e_y},\vect{e_z})$ of $\vect{\mathcal E}$ are 
\begin{align*}
 \{L_x,L_y\}&=-L_z\,,&\{L_y,L_z\}&=-L_x\,,& \{L_z,L_x\}&=-L_y\,;\\
 \{L_x,\varepsilon_x\}&=0\,,& \{L_x,\varepsilon_y\}&=-\varepsilon_z\,,& \{L_x,\varepsilon_z\}&=\varepsilon_y\,;\\
 \{L_y,\varepsilon_x\}&=\varepsilon_z\,,& \{L_y,\varepsilon_y\}&=0\,,& \{L_y,\varepsilon_z\}&=-\varepsilon_x\,;\\
 \{L_z,\varepsilon_x\}&=-\varepsilon_y\,,& \{L_z,\varepsilon_y\}&=\varepsilon_x\,,& \{L_z,\varepsilon_z\}&=0\,;\\
 \{\varepsilon_x,\varepsilon_y\}&=\frac{2E}{m^3k^2}L_z\,,&
 \{\varepsilon_y,\varepsilon_z\}&=\frac{2E}{m^3k^2}L_x\,,&  
 \{\varepsilon_z,\varepsilon_x\}&=\frac{2E}{m^3k^2}L_y\,.  
 \end{align*}
The family of fuctions on the phase space of the Kepler problem spanned by 
the components of $\vect L$ and of $\vect\varepsilon$ is not a Lie algebra, 
because the function $E$ appears in the right hand sides of equalities in the 
last line. Of course, if we replace $\vect\varepsilon$ by
$\displaystyle\frac{\vect\varepsilon}{\sqrt{-\zeta E}}$, we get a Lie algebra isomorphic to $so(4)$ when $E<0$
and to $so(3,1)$ when $E>0$. But in our opinion, this is not a good idea, because if we do this the energy level
$E=0$ is lost and the true geometric nature of the family of infinitesimal symmetries, 
on the whole phase space of the Kepler problem, is hidden. The above formulae make up the \emph{bracket table} of
smooth functions defined on the dual of a Lie algebroid (with zero anchor map), whose base is $\RR$ 
(spanned by the coordinate $E$) and standard fibre
$\RR^6=\RR^3\times\RR^3$ (spanned by the components of $\vect L$ and $\vect\varepsilon$ as coordinate functions).
 
\par\smallskip

The \emph{energy-momentum map} is the map $\mathcal J$, defined on the phase space of the Kepler problem, with values in
$\RR\times\vect{\mathcal E}\times\vect{\mathcal E}$, 
 $$(\vect r, \vect p)\mapsto{\mathcal J}(\vect r, \vect p)=\bigl(E(\vect r, \vect p), 
 \vect L(\vect r, \vect p), \vect\varepsilon(\vect r,\vect p)\bigr)\,.
 $$
The space $\RR\times\vect{\mathcal E}\times\vect{\mathcal E}$ will be called the \emph{energy-momentum space}. As seen above,
that space is the dual of a Lie algebroid; it should be considered 
as fibered over its first factor $\RR$, the fibre 
$\{e\}\times\vect{\mathcal E}\times\vect{\mathcal E}$ over each point $e\in\RR$ being equipped with a linear Poisson
structure which smoothly depends on $e$. With this structure, that fibre is the dual space of a Lie algebra which, 
as a vector space, can be identified with $\vect{\mathcal E}\times\vect{\mathcal E}$, with a Lie algebra bracket 
$[\ ,\ ]_e$ which smoothly depends on $e$. That bracket is easily deduced from the Poisson 
brackets of the components of $\vect L$ and $\vect\varepsilon$ indicated above: for each energy level $e$, the map
which associates to 
$(\vect{u_1},\vect{u_2})\in\left(\vect{\mathcal E}\times\vect{\mathcal E},[\ ,\ ]_e\right)$  
the Hamiltonian vector field with Hamiltonian
 $$(\vect r,\vect p)\mapsto \bigl\langle J(\vect r, \vect p), (\vect{u_1},\vect{ u_2})\bigr\rangle\,,$$
restricted to the energy level $e$ of the phase space,  
must be a Lie algebras homomorphism. This remark leads to the formula
 $$\bigl[(\vect{u_1},\vect{u_2}),(\vect{v_1},\vect{v_2}) \bigr]_e=
 \left(-\vect{u_1}\times\vect{v_1}+\frac{2e}{m^3k^2}\vect{u_2}\times\vect{v_2}\,,\,-\vect{u_1}\times\vect{v_2}
 +\vect{u_2}\times\vect{v_1}\right)\,.
 $$
By gluing together all the fibres $\vect{\mathcal E}\times\vect{\mathcal E}$ over all points $e\in \RR$, we get on
$\RR\times \vect{\mathcal E}\times\vect{\mathcal E}$ a Lie algebroid structure, whose dual is the energy-momentum space.
The above formula gives the bracket of two smooth sections $e\mapsto \bigl(\vect{u_1(e)},\vect{u_2(e)}\bigr)$ 
and $e\mapsto \bigl(\vect{v_1(e)},\vect{v_2(e)}\bigr)$ of that Lie algebroid. For each $e\in\RR$, 
$\bigl[(\vect{u_1},\vect{u_2}),(\vect{v_1},\vect{v_2}) \bigr](e)$
only depends on the values taken by $(\vect{u_1},\vect{u_2})$ and $(\vect{v_1},\vect{v_2})$ at point $e$, because 
the considered Lie algebroid has a zero anchor map.

\subsection{Action of a Lie algebroid on a symplectic manifold}\label{liealgebroidaction}
Let $\pi_A:A\to B$ be a Lie algebroid~\cite{cannasweinstein, mackenzie} whose base is a smooth manifold $B$
and whose anchor map is denoted by $\rho:A\to TB$. The bracket of two smooth sections $s_1:B\to A$ and $s_2:B\to A$
of $\pi_A$ is denoted by $\{s_1,s_2\}$. Let $\pi_M:M\to B$ be a surjective submersion of a smooth manifold $M$ onto the base
$B$. An \emph{action} of the Lie algebroid $\pi_A:A\to B$ on the fibered manifold $\pi_M:M\to B$ is a map
$s\mapsto X_s$ which associates, to each smooth section $s:B\to A$ of $\pi_A$, a vector field
$X_s$ on $M$, in such a way that for each smooth section $s$ of $\pi_A$ and each smooth function $f:B\to\RR$,
$$X_{fs}=(f\circ\pi_M)X_s\,,\quad T\pi_M\circ X_s=\rho\circ s\circ\pi_A\,,
 $$
and that for each pair $(s_1,s_2)$ of smooth sections of $\pi_A$,
 $$X_{s_1+s_2}=X_{s_1}+X_{s_2}\,,\quad[X_{s_1},X_{s_2}]=X_{\{s_1,s_2\}}\,.$$
Let us now assume that the manifold $M$ is endowed with a symplectic form $\omega$. The vector fields $X_s$ on $M$
associated to smooth sections $s$ of the Lie algebroid $\pi_A:A\to B$ cannot all be Hamiltonian, because if for some choice
of the section $s$, $X_s$ is Hamiltonian, the vector field $X_{fs}$, associated to the section $fs$, where $f:B\to\RR$ is a 
smooth function, will not in general be Hamiltonian. So one may wonder what should be a reasonable definition 
of a Lie algebroid Hamiltonian action on a symplectic manifold. The answer is suggested by the example of the Kepler problem.
\par\smallskip

In the last subsection, we described the Lie algebroid $\pi_A:A\to B$ of infinitesimal 
symmetries of the Kepler problem: $A=\RR\times\vect{\mathcal E}\times\vect{\mathcal E}$, $B=\RR$, $\pi_A:A\to B$
is the projection on the first factor, the anchor map is the zero map $A\to TB$. The vector bundle $\pi_A:A\to B$ being trivial,
each pair of vectors $(\vect{u_1},\vect{u_2})\in\vect{\mathcal E}\times\vect{\mathcal E}$ can be considered as a (constant) section
$s_{(\vect{u_1},\vect{u_2})}$ of our Lie algebroid. We defined the action of our Lie algebroid on the phase space of the Kepler problem 
by taking, as vector field $X_{s_{(\vect{u_1},\vect{u_2})}}$ associated to the section $s_{(\vect{u_1},\vect{u_2})}$, the Hamiltonian
vector field with Hamiltonian 
 $$(\vect r,\vect p)\mapsto\bigl\langle{\mathcal J}(\vect r,\vect p), s_{(\vect{u_1},\vect{u_2})}\bigr\rangle\,.$$
This condition unambiguously defines our Lie algebroid action (and, simultaneously, the bracket composition law of its sections) 
since the  module of its smooth sections is spanned by the constant sections.  
\par\smallskip

It would be interesting to see whether such a construction can be extended for more general Lie algebroids, with a nonzero anchor map.

\subsection{The $\mathbf S$ map}\label{smap}
The main disadvantage of Moser's regularization method is that it handles separately each energy level.
This disadvantage may be partially removed by the following procedure, which allows to handle together
all negative (resp., all positive) energy levels. However, negative, positive and zero energy levels still
cannot be handled together with that procedure.
\par\smallskip

Let us consider the quadric $Q_\rho$, defined by
 $$h^2+\zeta(x^2+y^2+z^2)=\zeta\rho^2\,,$$
where $\rho$ may take any positive value, and let $Q_R$ be the quadric defined by
 $$h^2+\zeta(x^2+y^2+z^2)=\zeta R^2\,,$$
where $R$ is a fixed positive quantity. To each point $M_R\in Q_R$, we associate the point $M_\rho\in Q_\rho$ such that
 $$\vect{OM_\rho}=\frac{\rho}{R}\,\vect{OM_R}\,.$$ 
We lift this diffeomorphism $Q_R\mapsto Q_\rho$ to the cotangent bundles, and we get a symplectic diffeomorphism
${\mathcal T}_\rho:T^*Q_R\to T^*Q_\rho$. This symplectic diffeomorphism associates to each pair $(M_R, \vect{W_R})$ made by a point $M_R\in Q_R$ and a vector $\vect{W_R}$ tangent to $Q_R$ at that point, the pair
$(M_\rho, \vect{W_\rho})$ made by a point $M_\rho\in Q_\rho$ and a vector $\vect{W_\rho}$ tangent to $Q_\rho$ at that point:
 $$\vect{OM_\rho}=\frac{\rho}{R}\,\vect{OM_R}\,,\quad \vect{W_\rho}=\frac{R}{\rho}\,\vect{W_R}\,.$$ 
For each $\rho>0$, we compose the symplectic diffeomorphism ${\mathcal T}_\rho: T^*Q_R\to T^*Q_\rho$
with the symplectic diffeomorphism $S_\rho:T^*(Q_\rho\backslash\{N_\rho\})\to T^*{\mathcal E}$ built in 
Subsection~\ref{cotangentlift}.
We obtain a family, indexed by $\rho>0$, of symplectic diffeomorphisms 
$S_{\rho,R}=S_\rho\circ{\mathcal T}_\rho:T^*(Q_R\backslash\{N_R\})\to T^*{\mathcal E}$: 
\begin{equation*}
\left\{
 \begin{aligned}
  \vect p&=\frac{\rho}{R-h_R}\,\vect {O\mu_R}\,,\quad\hbox{with}\quad h_R^2+\zeta\Vert \vect{O\mu_R}\Vert^2=R^2\,,\\
  \vect r&=\frac{R-h_R}{\rho}\,\vect {W_{3R}}+\frac{W_{hR}}{\rho}\,\vect{O\mu_R}\,,
 \end{aligned}
\right.
\end{equation*}
where we have set
 $$\vect{OM_R} =\vect{O\mu_R}+h_R\vect{e_h}\,,\quad \vect{W_R}=\vect{W_{3R}}+W_{hR}\vect{e_h}\,.$$ 
The formulae for the inverse transformation are
\begin{equation*}
\left\{
 \begin{aligned}
  \vect{O\mu_R}&=\frac{2R\rho}{\rho^2+\zeta p^2}\vect p\,,\\
  h_R&=R\,\frac{p^2-\zeta\rho^2}{p^2+\zeta\rho^2}\,,
     \quad\hbox{with}\quad p=\Vert\vect p\Vert=\sqrt{p_x^2+p_y^2+p_z^2}\,,\\
  \vect{W_{3R}}&=\frac{\rho^2+\zeta p^2}{2R\rho}\vect r-\frac{\zeta\vect r.\vect p}{R\rho}\vect p\,,\\
  W_{hR}&=\frac{\zeta\vect r.\vect p}{R}\,.
 \end{aligned}
\right.
\end{equation*}
The diffeomorphism $S_{\rho,R}^{-1}$ sends the subset of $T^*\bigl({\mathcal E}\backslash\{O\}\bigr)$ 
on which the energy is
$\displaystyle E(\vect r, \vect p)=-\frac{\zeta \rho^2}{2m}$ into the subset of $T^*Q_R$ on which
$\displaystyle\Vert\vect W_R\Vert=\frac{km^2}{R\rho}$.
Therefore, if $\rho_1$ and $\rho_2$ are two distinct
possible values of $\rho$, the images by $S_{\rho_1,R}^{-1}$ of the energy level
$\displaystyle E(\vect r, \vect p)=-\frac{\zeta \rho_1^2}{2m}$, and by
$S_{\rho_2,R}^{-1}$ of the energy level
$\displaystyle E(\vect r, \vect p)=-\frac{\zeta \rho_2^2}{2m}$, are disjoint.
By restricting each map $S_{\rho,R}^{-1}$ to the subset of $T^*\bigl({\mathcal E}\backslash\{O\}\bigr)$ 
on which the energy is
$\displaystyle E(\vect r, \vect p)=-\frac{\zeta \rho^2}{2m}$, and by gluing together these restricted maps for
all possible values of $\rho$, we obtain a unique diffeomorphism
$S^{-1}$ from the open subset of the cotangent bundle $T^*\bigl({\mathcal E}\backslash\{O\}\bigr)$ 
minus the zero section on which  the energy $E(\vect r, \vect p)$ is negative if $\zeta=1$, 
positive if $\zeta=-1$,
onto $T^*(Q_R\backslash\{N_R\})$ minus the zero section. This diffeomorphism, built by gluing together pieces 
of the symplectic diffeomorphisms $S_{\rho,R}^{-1}$
for all values of $\rho>0$, is no longer symplectic! It is given by the formulae, 
in which we no longer write the subscript $R$,
\begin{equation*}
\left\{
 \begin{aligned}
  \vect{O\mu}&=\zeta\,\frac{R\sqrt{\zeta r(2m^2k-rp^2)}}{m^2k}\,\vect p\,,\\
  h&=\frac{R(rp^2-m^2k)}{m^2k}\,,\\
  \vect{W_{3}}&=\zeta\,\frac{m^2k\vect r-r(\vect r.\vect p)\vect p}{R\sqrt{\zeta r(2m^2k-rp^2)}}\,,\\
  W_{h}&=\zeta\,\frac{\vect r.\vect p}{R}\,.
 \end{aligned}
\right.
\end{equation*}
We have set
 $$\vect{OM}=\vect{O\mu}+h\vect{e_h}\,,\quad \vect{W}=\vect{W_3}+W_h\vect{e_h}\,.$$
 The transformed Hamiltonian $(S^{-1})^*E=E\circ S$ is
$$H=E\circ S=-\frac{\zeta k^2m^3}{2R^2}\,\frac{1}{\Vert\vect W\Vert^2}\,.$$
Up to the constant factor $\displaystyle\frac{\zeta k^2m^3}{R^2}$, it is the \emph{Delaunay Hamiltonian}~\cite{cushmanbates,cushmanduistermaat},
\emph{i.e.}, the Kepler Hamiltonian in Delaunay coordinates. This property is related to the fact that Delaunay variables are 
action-angle variables for the Kepler problem.
\par\smallskip
In the rest of the paper we write $Q$ for $Q_R$. The north pole of $Q$ is denoted by $N$ instead of $N_R$.
\par\smallskip
 
Let $X_E$ be the Hamiltonian vector field on $T^*({\mathcal E}\backslash\{0\})$ with Hamiltonian $E$, 
\emph{i.e.} the Hamiltonian vector field of the Kepler problem, and let $X_H$ be the Hamiltonian
vector field, on $T^*\bigl(Q\backslash\{N\}\bigr)$ minus the zero section, whose Hamiltonian 
$H=E\circ S$ is given by the above formula. The direct image $(S^{-1})_*(X_E)$ of the vector field 
$X_E$ by the diffeomorphism $S^{-1}$ is not equal to $X_H$, since $S^{-1}$ is not symplectic. 
A short calculation leads to its expression,
 $$(S^{-1})_*(X_E)= g\,X_H\,,\quad\text{with}\quad g=\frac{R}{h-R}\,.$$
The function $g$ is smooth on $T^*\bigl(Q\backslash\{N\}\bigr)$ minus the zero section, and becomes singular
when $h=R$, \emph{i.e.} on the fibre over the north pole $N$. 
We see that the map $S^{-1}$ sends the Hamiltonian vector field of the Kepler problem 
to a conformally Hamiltonian vector field, with $H$ as Hamiltonian and with $g$ as conformal factor.

\subsection{The flow of $\mathbf {X_H}$}\label{flowofXH}
Unlike the conformal factor $g$, which becomes singular 
on the cotangent space to $Q$ at the north pole, the Hamiltonian $H$ is smoothly
defined on the whole cotangent bundle $T^*Q$ minus the zero section. 
The associated Hamiltonian vector field $X_H$ is complete. 
Its flow, defined on $\RR\times \bigl(T^*Q\backslash\{\text{zero section}\}\bigr)$,
 $$\Phi_{X_H}\left(s,\bigl(\vect{OM(0)},\vect{W(0)}\bigr)\right)=\bigl(\vect{OM(s)},\vect{W(s)}\bigr)
 $$
has slightly different expressions for $\zeta=1$ (negative energy)
and for $\zeta=-1$ (positive energy). Observe that $\Vert\vect W\Vert$ is a first integral of $X_H$.
 \begin{equation*}
 \text{For\ }\zeta=1\ (E<0)\,,\ 
 \left\{
 \begin{aligned}
 \vect{OM(s)}&=\cos(\lambda s)\,\vect{OM(0)}
               +\frac{R}{\Vert\vect{W(0)}\Vert}\,\sin(\lambda s)\,\vect{W(0)}\,,
               \\
 \vect{W(s)}&=-\frac{\Vert\vect{W(0)}\Vert}{R}\,\sin(\lambda s)\,\vect{OM(0)}+\cos(\lambda s)\,\vect{W(0)}\,.
 \end{aligned}
 \right.
 \end{equation*}
 \begin{equation*}
 \text{For\ }\zeta=-1\ (E>0)\,,\ 
 \left\{
 \begin{aligned}
 \vect{OM(s)}&=\cosh(\lambda s)\,\vect{OM(0)}
               -\frac{R}{\Vert\vect{W(0)}\Vert}\,\sinh(-\lambda s)\,\vect{W(0)}\,,
               \\
 \vect{W(s)}&=-\frac{\Vert\vect{W(0)}\Vert}{R}\,\sinh(-\lambda s)\,\vect{OM(0)}+\cosh(\lambda s)\,\vect{W(0)}\,.
 \end{aligned}
 \right.
 \end{equation*}
 We have set $\displaystyle\lambda=\frac{k^2m^3}{R^3\Vert\vect{W(0)}\Vert^3}$.

\subsection{A symplectic diffeomorphism}\label{symplecticdiff}
At the end of Subsection~\ref{smap}, we have seen that the image by $S^{-1}$ of the Hamiltonian 
vector field $X_E$ of the Kepler problem is a conformally Hamiltonian vector field 
$gX_H$ defined on an open dense subset of $T^*Q\backslash\{\text{zero section}\}$,
whose Hamiltonian $H$ is smoothly defined on the whole $T^*Q\backslash\{\text{zero section}\}$; 
the corresponding Hamiltonian vector field $X_H$ is complete.
Moreover the conformally Hamiltonian vector field $gX_H$ is also Hamiltonian, with Hamiltonian $H$,
not for the canonical symplectic form of $T^*Q$, but for the pull-back by $S$ of the canonical 
symplectic form of $T^*{\mathcal E}$. We have seen in Subsection~\ref{levi-civitaparameter} that for every
solution $t\mapsto\bigl(\vect{r(t)},\vect{p(t)}\bigr)$ of the equations of motion of the Kepler problem,
 $$\frac{d}{dt}\left(\frac{\vect{p(t)}.\vect{r(t)}-2E\bigl(\vect{r(t)},\vect{p(t)}\bigr)t}{mk}\right)
 =\frac{1}{r(t)}\,.$$ 
The conformal factor $g$, whose expression is given at the end of Subsection~\ref{smap},
may also be expressed as
 $$g=\frac{mk}{2}\,S^*\left(\frac{1}{rE(\vect r, \vect p)}\right)\,.$$
Therefore
 $$\frac{d}{dt}\left(\frac{\vect{p(t)}.\vect{r(t)}-2E\bigl(\vect{r(t)},\vect{p(t)}\bigr)t}
 {2E\bigl(\vect{r(t)},\vect{p(t)}\bigr)}\right)=g\circ S^{-1}\,.
 $$
All the integral curves of the conformally Hamiltonian vector field $gX_H$ can be written
 $$t\mapsto S^{-1}\circ\bigl(\vect{r(t)},\vect{p(t)}\bigr)\,,$$
where $t\mapsto\bigl(\vect{r(t)},\vect{p(t)}\bigr)$ is a solution of the equations of motion of the Kepler problem.
We see that the pull-back by $(\id_\RR,S)$ of the function, defined on the product with $\RR$ 
of the open subset of the phase space of the Kepler problem on which the energy $E$ is negative if
$\zeta=1$, positive if $\zeta=-1$,
  $$\frac{\vect p\cdot\vect r-2E(\vect r, \vect p)\,t}{2E(\vect r, \vect p)}$$
has the properties of the function $\sigma$ of Theorem~\ref{result1} of Subsection~\ref{confhamfields}.
That Theorem can therefore be applied. It proves that by composing the symplectic diffeomorphism
$S^{-1}$ with the flow $\Phi_{X_H}$ of the Hamiltonian vector field $X_H$, for suitably chosen values of the independent variable $s$, we get a symplectic diffeomorphism from the phase space of the Kepler problem, restricted either to negative, or to positive energies,  onto an open subset of $T^*Q$. The explicit expression 
of this symplectic diffeomorphism is
 $$\left(\vect p,\vect r\right)\mapsto\left(\vect{OM(s)},\vect{W(s)}\right)\,,\quad\text{with}
 \quad s=\frac{-\vect p\cdot\vect r}{2E\left(\vect r,\vect p\right)}\,.$$
The quantities $\bigl(\vect{OM(s)},\vect{W(s)}\bigr)$ are given, as functions of
$\bigl(\vect{OM(0)},\vect{W(0)}\bigr)$, by the formulae at the end of Subsection~\ref{flowofXH}.
We have to set
 $$\vect{OM(0)}=\vect{O\mu}+h\vect{e_h}\,,\quad \vect{W(0)}=\vect{W_3}+W_h\vect{e_h}\,.$$
The formulae at the end of Subsection~\ref{smap} give the expressions of
$\vect{O\mu}$, $h$, $\vect{W_3}$ and $W_h$ as functions of $\vect{r}$ and $\vect{p}$.  
Finally we see that the symplectic diffeomorphism so obtained is that of Gy\"orgyi, Ligon and Schaaf.
R.~Cushman and L.~Bates, in chapter II of the book~\cite{cushmanbates}, offer a detailed discussion of its properties,
notably its behaviour with respect to the Liouville $1$-forms of the cotangent spaces of the Kepler 
configuration manifold and of the $3$-dimensional sphere (or hyperboloid).

\section{Conclusion and perspectives}
The properties of the diffeomorphism from the phase space of the Kepler problem
to the cotangent bundle of a sphere or of an hyperboloid are explained, in Subsection~\ref{smap}, 
from a point of view other than that used by Cushman and Duistermaat~\cite{cushmanduistermaat}. 
Our explanation is founded on a very natural property of conformally 
Hamiltonian vector fields. However, the procedure we used, as well as those used by Gy\"orgyi~\cite{gyorgyi},
Ligon and Schaaf~\cite{ligonschaaf} remains unsatisfactory in the fact that it handles separately 
negative and positive energy levels. The space of motions of the Kepler problem is connected, 
since by varying slowly the energy level, an elliptic motion can be transformed into a parabolic, 
then into an hyperbolic motion. By completely different methods, J.-M.~Souriau~\cite{souriautorino}
proposed a global description of the manifold of motions of the Kepler problem, of its regularization and
of its global and infinitesimal symmetries. It should be possible to do something similar at the level
of the phase space instead of at the level of the manifold of motions.
\par\smallskip

The fact that the Kepler vector field is Hamiltonian with respect to a symplectic form and conformally Hamiltonian 
with respect to another syplectic form is remarkable; it is probably related to the complete integrability of that vector field, as suggested by 
an anonymous referee. In~\cite{MaciePryTsi}, A.~Maciejewski, M.~Prybylska and A.~Tsiganov have used conformally Hamiltonian vector fields
within the theory of bi-Hamiltonian systems, to build completely integrable systems. 
\par\smallskip

We have shown that the infinitesimal symmetries of the Kepler problem form a Lie algebroid rather than a Lie algebra. Actions of Lie algebroids on
symplectic manifolds is a subject which seems to us interesting; it should be nice to have other examples of Hamiltonian vector fields
whose infinitesimal symmetries form a Lie algebroid, maybe with a non-zero anchor map.
\par\smallskip

A.~Douady and M.~Lazard~\cite{douadylazard} have shown that Lie algebroids with a zero anchor map are integrable into Lie groupoids; it would be interesting to write down explicitly the action of a Lie groupoid with the groups ${\rm SO}(4)$, ${\rm SE}(4)$ and ${\rm SO}(3,1)$ as isotropy groups
on a regularized phase space of the Kepler problem.

\section*{Acknowledgments}
I am much indebted towards my colleagues and friends of the \lq\lq Centre International de Th\'eories
Variationnelles\rq\rq, the International Workshop on Differential Geometric Methods in Mechanics
and the Seminar of Hamiltonian Geometry for helpful discussions, towards Alan Weinstein who, 
after listening to a presentation of a preliminary version of this work, made several enlighting observations,  
and towards Alain Guichardet and Alain Albouy who communicated me their unpublished papers on the Kepler problem.
I address my thanks to the two anonymous referees whose judicious observations and suggestions were very helpful.  


\medskip
Received xxxx 20xx; revised xxxx 20xx.
\medskip

\end{document}